\documentclass{article}
\usepackage[latin1]{inputenc}
\usepackage[T1]{fontenc}
\usepackage{amssymb}
\usepackage{showlabels}
\usepackage{longtable}

\usepackage{color}
\definecolor{DarkOlive}{rgb}{0.1047,0.2412,0.0064}
\definecolor{FireBrick}{rgb}{0.5812,0.0074,0.0083}
\definecolor{RoyalBlue}{rgb}{0.0236,0.0894,0.6179}
\definecolor{RoyalGreen}{rgb}{0.0236,0.6179,0.0894}
\definecolor{RoyalRed}{rgb}{0.6179,0.0236,0.0894}
\definecolor{LightBlue}{rgb}{0.8544,0.9511,1.0000}
\definecolor{Black}{rgb}{0.0,0.0,0.0}
\definecolor{MidnightBlue}{rgb}{0.0035,0.0020,0.1363}
\definecolor{FireBrick3}{rgb}{0.5812,0.0074,0.0083}
\definecolor{FireBrick4}{rgb}{0.2156,0.0023,0.0035}
\definecolor{Blue2}{rgb}{0.0000,0.0000,0.8644}
\definecolor{Navy}{rgb}{0.0000,0.0000,0.1927}
\definecolor{MediumBlue}{rgb}{0.0000,0.0000,0.6179}
\usepackage[hyphens]{url} 
\usepackage[
        a4paper=true,bookmarks=false,pdftitle={Generic Computations in
        Finite Groups of Lie Type},
        colorlinks=true,backref=false,breaklinks=true,linkcolor=MediumBlue,
        citecolor=FireBrick3,filecolor=RoyalRed,
        urlcolor=Blue2,pagecolor=MediumBlue]{hyperref}
\usepackage{tikz}

\newcommand{\F}{\mathbb{F}}

\newcommand{\GG}{\ensuremath{{}^2\!E_6(2)}}
\newcommand{\bG}{\ensuremath{{\boldmath \textbf{G}}}}
\newcommand{\ATLAS}{\textsf{ATLAS}}
\newcommand{\CHEVIE}{\textsf{CHEVIE}}
\newcommand{\GAP}{\textsf{GAP}}
\newcommand{\Magma}{\textsf{Magma}}
\newcommand{\CTblLib}{\textsf{CTblLib}}
\newcommand{\Irr}{{\textrm{Irr}}}

\begin{document}
\title{Characters and Brauer trees of the covering group of \GG}
\author{Frank Lübeck}
\date{}
\maketitle
\begin{abstract}
Let $G$  be the finite simple  Chevalley group of  type \GG. It has  a Schur
multiplier of type  $C_2^2 \times C_3$. We determine  the ordinary character
tables of the central extensions $3.G$,  $6.G$, $(2^2\times 3).G$ of $G$ and
their extensions by  an automorphism of order $2$, that  is $3.G.2$, $6.G.2$
and $(2^2\times 3).G.2$.

Furthermore we  determine all  Brauer trees  of all  groups of  type $Z.G.A$
(where $Z$ is central in $Z.G \lhd Z.G.A$ and $A \cong Z.G.A/Z.G$) for which
the ordinary character table is known.

\end{abstract}

\section{Introduction}

Let $p$  be a  prime and  \bG\ a  simple reductive  algebraic group  over an
algebraic closure of the finite field $\F_p$ with $p$ elements. Let $q$ be a
power of $p$ and  \bG\ be defined over $\F_q$. We  denote $\bG(q)$ the group
of $\F_q$-rational points--a (twisted or untwisted) finite group of Lie type. 
Let $Z$ be the center of $\bG(q)$.

If \bG\  is of  simply connected  type then except  for finitely  many cases
$\bG(q)/Z$ is  a finite simple group.  If $\bG(q)/Z$ is simple  then in most
cases $Z$ is  isomorphic to the Schur multiplier of  the simple quotient and
$\bG(q)$ is a  universal covering. There are $18$ cases  with an exceptional
Schur multiplier,  see~\cite[3. and Table~5.]{CC85}.  For only one  of these
cases the  (ordinary) character  table of  the full  covering group  was not
determined by the  \ATLAS\ project, this is the case  where $\bG(q)/Z$ is of
type \GG.  In the first  part of this paper  we describe the  computation of
this missing character  table. This also yields the character  table of this
group  extended by  an outer  automorphism of  order $2$  and of  some other
related groups.

Now let $G$ be the finite simple group of type \GG. 

\textbf{Previously  known character  tables.}  The character  tables of  the
central extension  $2^2.G$ of  $G$ with  the exceptional  part of  the Schur
multiplier  (and so  of its  factor groups  $2.G$ and  $G$), as  well as  of
the  larger  groups  of  form $G.2$,  $G.3$,  $G.S_3$,  $2.G.2$,  $2^2.G.2$,
$2^2.G.3$  and $2^2.G.S_3$  were  available in  the  \GAP{} character  table
library~\cite{CTblLib1.1.3}, the  table of  $2.G.2$ is  also printed  in the
\ATLAS~\cite[p.191]{CC85}.

The  essential  step to  get  the  character  table  of the  covering  group
$(2^2\times 3).G$ of $G$ is to find the table of its quotient of type $3.G$.
The group $3.G$ arises as a twisted  finite group of Lie type $\bG(2)$ where
\bG\ is  the simple simply connected  reductive group of type  $E_6$ over an
algebraic closure of the field with $2$  elements. For this group we can use
Deligne-Lusztig  theory  (see~\cite{DL76},  \cite{Ca85}) to  construct  some
faithful  characters, which  together with  the  known table  of the  simple
quotient $G$ enables us to compute the  whole table of $3.G$. We explain the
details of the construction in section~\ref{Constr3G}.

As an  application we describe  in section~\ref{trees} the  determination of
all Brauer  trees (that is  the decomposition  numbers for all  prime blocks
with  a  cyclic  defect  group)  of  $(2^2\times  3).G$  and  the  available
extensions of this and its factor groups by outer automorphisms.

Since the newly computed character  tables are huge (e.g., $(2^2\times 3).G$
has $934$ conjugacy classes and  irreducible characters) we will not include
the  new  character  tables  in  this  article.  The  tables  are  available
in  newer   versions  of  the  \GAP-package   \CTblLib,  see~\cite{GAP4.4.9}
and~\cite[since  version~1.2.0]{CTblLib1.1.3},  and  were already  used  for
various applications. Information on conjugacy classes and character degrees
of $3.G$ is also given on the webpage~\cite{web2e62}.

We note that  the \ATLAS-tables which are relevant here  could be recomputed
and checked recently~\cite{Rev16} using  some explicit representations and a
generic function in \Magma~\cite{MAGMA}. But the case of the table for $3.G$
seems  out  of reach  of  this  generic  function currently.  Therefore,  we
document the  computation of  the tables of  groups with  simple composition
factor $G$ which are not contained in the \ATLAS\ in this note.

The character tables  of the related groups  $3.G.3$, $3.G.S_3$, $(2^2\times
3).G.3$ and  $(2^2\times 3).G.S_3$  are still  not known.  Their computation
will be addressed elsewhere.

\section{The character table of {\boldmath $3.G$}}\label{Constr3G}

As  mentioned in  the Introduction  the group  $3.G$ can  be constructed  as
twisted  group $\bG(2)$  of  $\F_2$-rational points  of  a simply  connected
algebraic group $\bG$ of type $E_6$ over an algebraic closure of $\F_2$.

To  compute  its character  table  we  use  the  following data  which  were
previously known:
\begin{itemize}
\item[(a)]  The  character  table  of  the  simple  quotient  $G$  from  the
{\ATLAS}~\cite{CC85}.  This  contains  the   centralizer  orders  (or  class
lengths),  the  power maps  (that  is  for  each  conjugacy class  $C$  with
representative $g$ and  integer $n$ one can read off  the conjugacy class of
$g^n$, in particular the orders  of representatives of the conjugacy classes
are known) and the character values.
\item[(b)] A parameterization  of the conjugacy classes of  $\bG(q)$ for all
prime powers $q$ with $q \equiv 2 \textrm{ mod }3$ and part of the values of
unipotent  characters  were determined  and  used  in~\cite{HL97}. This  can
easily be  specialized to the $346$  conjugacy classes in case  $q=2$. Those
computations also  gave explicit representatives for  the semisimple classes
(those  with elements  of odd  order)  in some  fixed maximal  torus of  the
algebraic  group $\bG$,  in particular  we get  the element  orders for  the
semisimple classes.
\item[(c)] A list  of all irreducible character degrees of  $\bG(q)$ for all
prime powers $q$ was computed for  the results in~\cite{L01}, again this can
easily be specialized to $q=2$.
\end{itemize}
The specialized data are available on~\cite{web2e62}.

To find some values of faithful characters for the group $3.G$ we extend~(b)
and compute  the Deligne-Lusztig characters of  $\bG(2)$. For this we  use a
practical variant of the  character formula~\cite[7.2.8]{Ca85}, as described
in~\cite[Satz~2.1(b)]{L93}.  In principle  this  allows one  to compute  the
generic values of  Deligne-Lusztig characters, that is  a parameterized form
of the values of $\bG(q)$ for all prime powers $q \equiv 2 \textrm{ mod }6$.
But this would need too much  memory and computation time, even with current
computers. Fortunately, only a small part of this generic character table is
relevant for  $q=2$ (for example,  only a  few conjugacy classes  of maximal
tori  of  $\bG(2)$ contain  regular  elements;  more  precisely 342  of  the
494  columns of  the  generic character  table are  not  needed for  $q=2$).
Computing just the  part of the generic table  of Deligne-Lusztig characters
relevant for $q=2$  is feasible. The values of the  Green functions occuring
in  the  character  formula  are  needed for  groups  of  type  ${}^2\!E_6$,
$D_4$  and  $A_l$,  $l  \leq  5$. They  are  available  from~\cite{M93}  and
{\CHEVIE}~\cite{GH96}.

Some  irreducible  characters  of  $\bG(2)$ are  known  linear  combinations
of   Deligne-Lusztig  characters,   this  uses   Lusztig's  parameterization
of   irreducible  characters~\cite[4.22]{Lu85}.   More  details   are  given
in~\cite[7.]{L93}.  There are  also some  further uniform  almost characters
which are  known as explicit  linear combinations of  irreducible characters
and of  Deligne-Lusztig characters;  the multiplicities of  irreducibles are
rational numbers,  and we can multiply  by some integers to  get generalized
characters.  This way we get:

\begin{itemize}
\item[(d)] The  values of  about 200 irreducible  characters of  $\bG(2)$ in
terms of the parameterization in~(b), some of them are trivial on the center
(that is they lead to known characters of the simple quotient $G$), and some
are not (that is they are faithful characters of $3.G$).
\item[(e)] The values of about 150  generalized characters of $\bG(2)$.
\end{itemize}

It remains to  explain how to merge the information  from~(a) to~(e) to find
the whole  character table of  $\bG(2)$. This involves  further computations
with \GAP~\cite{GAP4.4.9}  which provides  powerful functions  for computing
with character tables.

First we need the fusion of the conjugacy classes as described in~(b) to the
classes of  the \ATLAS-table  in~(a). There are  two possible  behaviours of
a  class,  either  an  element  multiplied by  the  center  elements  yields
representatives of  one or  of three  different classes;  that is  under the
canonical map $3.G \to G$ the preimage  of each class of $G$ contains either
one or  three classes. Using  the element  orders and centralizer  orders in
both cases  and identifying some  irreducible characters from~(d)  which are
trivial on  the center  with characters in~(a)  only very  few possibilities
remain. Namely,  on a few  tuples of  classes all characters  in~(d) and~(e)
have the same value, in these cases we just choose some fusion.

Now we  can lift all  irreducible characters of $G$  as given in~(a)  to the
table of $3.G$.

It turns out that the class functions we have found so far are sufficient to
determine all power maps for $3.G$ using the \GAP-function for computing all
power maps  compatible with  a set  of given  characters (there  are several
possibilities but these  are all equivalent modulo  simultaneous renaming of
some conjugacy classes and irreducible characters).

The next step is to compute many tensor products of known virtual characters
and  to apply  \GAP  s implementation  of the  LLL-algorithm  to find  class
functions of  small norm in  the lattice  of all available  class functions.
With this technique we find 322 of the 346 irreducible characters as well as
20 virtual characters of norm 2 or 3 and 4 virtual characters of norm 24.

There are some further standard tricks  to produce more class functions from
known ones, like computing symmetrizations  or inducing class functions from
subgroups (e.g., in  our case there is a subgroup  of type $3.Fi_{22}$). But
all of these do not improve the state as described above.

To  find the  remaining irreducible  characters  we can  now take  advantage
of~(c),  the list  of known  character degrees.  This tells  us that  we are
missing irreducible  characters of  degrees 7194825,  1929727800, 4583103525
and 11972188800,  each occuring 6  times. It turns out  that for all  of our
virtual characters  of norm 2  or 3 there is  only one possibility  to write
their degree as  sum or difference of  the missing degrees. So,  we know the
degrees  and the  multiplicities ($\pm  1$)  of their  constituents. Now  we
consider  the  scalar products  of  our  non-irreducible virtual  characters
with  themselves and  their complex  conjugates. This  shows that  all their
constituents are different  from their complex conjugates, and  that any two
of the virtual characters of norm 2  and 3 and their complex conjugates have
at  most one  common constituent.  This yields  a labelling  of the  missing
characters and a decomposition of the 20  virtual characters of norm 2 and 3
in terms  of this labelling.  For the  remaining four virtual  characters of
norm 24 there are several  possible decompositions which are compatible with
the  computed scalar  products. But  the  number of  possibilities is  small
enough that  we can try  all of  them and compute  the potential set  of all
irreducible  characters. Some  of  these possibilities  could  be ruled  out
easily,  because  some  random  tensor  product  has  scalar  products  with
irreducibles which are not non-negative integers.

In the end there  are 12 possibilities left, and it turns  out that they are
all  equivalent modulo  table automorphism  (that is,  the resulting  tables
are  the same  modulo  simultaneous permutations  of  conjugacy classes  and
irreducible characters).

\subsection{\boldmath Character tables of $6.G$, $(2^2\times 3).G$,
$3.G.2$, $6.G.2$ and $(2^2\times 3).G.2$}

Having constructed the character table of $3.G$ as described above and using
the previously  known character tables  for $2.G$, $2^2.G$,  $G.2$, $2.G.2$,
$2^2.G.2$  it turns  out to  be straightforward  to compute  the new  tables
mentioned in  the header. This can  be done with utility  functions from the
character table \GAP-package \CTblLib~\cite{CTblLib1.1.3}.

Using   the   known   tables   of    $G$,   $2.G$   and   $3.G$   as   input
we   can   now   compute   the    table   of   $6.G$   with   the   function
\texttt{CharacterTableOfCommonCentralExtension}  (this function  computes as
many irreducibles of the  common extension as it can find,  in this case all
irreducibles of $6.G$ are found).

In the next step  we can construct the table of  the group $(2^2\times 3).G$
with the  function \texttt{PossibleCharacterTablesOfTypeV4G}. The  input for
this function are the  tables of $3.G$ and $6.G$ as  well as the permutation
induced by an  automorphism of order $3$ on the  conjugacy classes of $3.G$.
Using the  function \texttt{AutomorphismsOfTable} on  the table of  $3.G$ it
turns out that the  result contains a unique subgroup of  order $3$ and this
way  we find  the needed  permutation. It  takes several  hours to  find the
possible  tables for  the larger  extension, it  turns out  that there  is a
unique possibility which then must be the table we are looking for.

Finally,  the  tables  for  the  groups of  form  $M.G.2$  can  be  computed
with the  function \texttt{PossibleCharacterTablesOfTypeMGA}.  This function
gets  as  input  the  tables  of  $M.G$,   $G$  and  $G.2$  and  a  list  of
orbits  of   the  group  $M.G.2$  acting   on  the  classes  of   $M.G$.  In
all  of  our  cases  these  orbits   are  easily  found  with  the  function
\texttt{PossibleActionsForTypeMGA} which only  returns a unique possibility.
Applying  \texttt{PossibleCharacterTablesOfTypeMGA} on  this input  yields a
unique possible result and so this must be the correct table.

\section{Determination of Brauer trees}\label{trees}

As an  application of the new  character tables constructed in  the previous
section we  want to determine all  Brauer trees of the  corresponding groups
which encode for primes $l$  the $l$-modular decomposition numbers of blocks
with non-trivial cyclic defect group.

For  the   following  facts  about   Brauer  trees  and   their  computation
from   character  tables   we  refer   to  the   first  three   chapters  of
the  book~\cite{BBB}.   By  results  of  Feit~\cite{F84}   and  K\"ulshammer
(unpublished,  see~\cite[1.3]{BBB}) the  shape  of any  Brauer  tree of  any
finite group occurs as Brauer tree of a central extension of an automorphism
group of  a finite simple  group. The latter are  now known for  most finite
simple groups, the cases considered here being among the few missing ones.

We briefly recall some  basic facts about Brauer trees. Let  $H$ be a finite
group  and $l$  be a  prime, and  let $(K,R,k)$  be a  splitting $l$-modular
system for $H$.  The set of irreducible characters  $\Irr(H)$ is partitioned
into $l$-blocks,  let $B = \{\chi_1,  \ldots, \chi_k\}$ be such  a block. We
assume that the  defect group of $B$ (a subgroup  of some Sylow $l$-subgroup
of $H$) is a non-trivial cyclic group  of order $l^d$. Let $e$ be the number
of irreducible $l$-modular Brauer characters of $B$. Then we have:

\begin{itemize}
\item[(a)] $e \mid (l-1)$.
\item[(b)] There may be several characters  in $B$ with the same restriction
to $l$-regular  classes. If this  happens, the corresponding  characters are
called  the exceptional  characters of  $B$ and  there are  $m =  (l^d-1)/e$
exceptional  and  $e$  non-exceptional   characters  in  $B$  (so  $k=e+m$).
Otherwise, there are no exceptional characters and $k = e+1$ non-exceptional
characters.
\item[(c)] All projective indecomposable characters corresponding to $B$ are
of the form  $\chi + \chi'$ where  $0\neq \chi \neq \chi'\neq  0$ and $\chi$
and $\chi'$  are either the  sum of all exceptional  characters in $B$  or a
non-exceptional character in $B$.
\item[(d)] The  decomposition matrix of $B$  can be encoded in  a graph. Its
$e+1$ vertices are labelled by the non-exceptional characters in $B$ and, if
there are  exceptional characters,  the sum  of all  exceptional characters.
There is  an edge  joining two vertices  $\chi$ and $\chi'$  if and  only if
$\chi + \chi'$ is projective.
\item[(e)] The graph defined  in (d) is a tree, that is  it is connected and
has $e$ edges. This is called the Brauer tree of $B$.
\item[(f)]  Let $\alpha$  be  an  automorphism of  $H$.  We  also denote  by
$\alpha$ the  induced map  on class  functions of $H$,  and for  a character
$\chi$  we write  $\bar\chi$ for  its complex  conjugate. Assume  that $\chi
\mapsto \bar\chi^\alpha$  maps $B$  to $B$.  Then this  map induces  a graph
automorphism on  the Brauer tree  of $B$ and  the subgraph of  the invariant
vertices ($\chi = \bar\chi^\alpha$) forms a line (in case $\alpha=1$ this is
called the real stem of the tree).
\item[(g)]  Let $\alpha$  be an  automorphism as  in~(f) which  maps $B$  to
another  block $B'$.  Then $\alpha$  induces  a graph  isomorphism from  the
Brauer tree of $B$ to the Brauer tree of $B'$.
\end{itemize}

Now we describe  the strategy which allowed  us to find the  Brauer trees we
are considering  here. Let  $H$ be a  group such that  we know  its ordinary
character table (the character values,  the centralizer orders and the power
maps)  and let  $l$ be  a  prime divisor  of $H$  (otherwise all  $l$-blocks
contain  a  single  irreducible  character).  We assume  that  we  have  the
character table available  in \GAP~\cite{GAP4.4.9} such that we  can use the
\GAP-functions for computing with character tables.

(1)  The  $l$-blocks  of the  character  table  of  $H$  and the  orders  of
the  corresponding defect  groups  can be  computed  with the  \GAP-function
\texttt{PrimeBlocks}. (This uses that two irreducibles are in the same block
if and only if  their central characters modulo $l$ are  equal, the order of
the defect  group can be  found from the  character degrees.) If  the defect
group is  of order  $l$ then it  is clearly cyclic  and non-trivial.  In the
cases we consider here it is easy to see that all other blocks do not have a
non-trivial  cyclic defect  group (because  it is  trivial or  the numerical
conditions in facts~(a), (b) are not fulfilled).

(2) If $G$  is the finite simple group  of type $\GG$ and $H$ is  one of the
groups $G$, $2^2.G$ or $(2^2\times 3).G$  then $H$ has an outer automorphism
group of type  $S_3$~\cite{CC85}. In \GAP\ we can  compute the automorphisms
of the character table of $H$ (a permutation of the classes, compatible with
power maps, that leaves  the table invariant). It turns out  that this has a
unique  subgroup of  type $S_3$.  Since  the outer  automorphism group  acts
faithfully on the conjugacy classes we find the explicit action of the outer
automorphisms of  $H$ on  $\Irr(H)$, as  well as the  induced action  on the
$l$-blocks of $H$.

The following  steps are applied  to one block $B$  of cylic defect  in each
orbit under  the outer automorphisms of  $H$ (this is sufficient  because of
fact~(g)).

(3) We start with listing for each  vertex in the tree all possible vertices
they may  be connected to.  For the initial list  we take into  account that
projective characters are zero on $l$-singular classes. So, if one character
has a non-zero value on some $l$-singular  class it can only be connected to
other characters  which have the negative  of that value on  the same class.
(If the tree has the shape of a star, that is, there is one vertex connected
to all the others, then we have already found the tree in this step.)

(4) Whenever we  find a new edge  of the tree during the  following steps we
may  be  able to  reduce  the  possibilities  for  further edges  using  the
facts~(e)  (there cannot  be a  further edge  between vertices  in the  same
connected component with  respect to the known edges)  and~(f) (an invariant
vertex can  only be connected to  at most two other  invariant vertices). If
the number of known edges and possible further edges is $e$, we are done.

In the situation of fact~(f) each  edge between $\chi$ and $\chi'$ involving
a non-invariant  character implies a further  edge between $\bar\chi^\alpha$
and $(\bar{\chi'})^\alpha$.

(5) We use that the defect-zero  characters (the single characters in blocks
with  trivial defect  group)  are  projective and  that  tensor products  of
projective characters  with arbitrary characters are  again projective. This
way we can easily compute a huge number of projective characters and compute
their scalar products with the irreducibles  (and the sum of the exceptional
characters) in  our block  $B$. For  each projective this  yields a  list of
multiplicities $m_i$  for each vertex  $v_i$ in our  tree ($m_i$ is  the sum
of  the  multiplicties of  all  projective  indecomposable characters  which
correspond to an  edge of the tree involving $v_i$).  Let $\{v_j\mid\; j \in
J\}$ the subset  of vertices which are possibly connected  to $v_i$, and let
$j' \in J$ such  that $m_{j'}$ is maximal among $\{m_j\mid\;  j \in J\}$. If
now $\sum_{j\in  J, j \neq j'}  m_j < m_i$  then we can conclude  that $v_i$
must  be connected  to $v_{j'}$;  we  have found  an  edge of  the tree.  In
particular we  have found a  new projective character corresponding  to this
edge, we  use this for further  iterations of this step.  (Heuristically, we
find a few (and  sometimes all) edges very quickly in  this step but nothing
new later. So, we stop this proceduce when we have not found new edges for a
while.)

(6)  In all  our cases  we find  enough  edges in~(5)  such that  it is  now
feasible to enumerate all trees which  are consistent with the edges already
found. Many of the  possible trees can be ruled out  easily by computing the
degrees of the irreducible Brauer characters  which are implied by the tree;
these  must be  positive  integers  but incorrect  trees  often yield  other
numbers.

(7) For each possible tree left in step~(6) we now compute again some random
projective characters (by tensoring known projectives with irreducibles) and
check  if the  multiplicities  with  the characters  in  $B$ are  consistent
with  the tree  (the multiplicities  of the  projective indecomposables  are
recursively determined by  the tree, starting from the leaves  of the tree).
This quickly rules out more trees.

(8) In very few cases we need to induce projective indecomposable characters
from a subgroup  of index~2 or~3 to find  a tree. In one case  this was also
not enough, but the induced projective characters with only four irreducible
constituents allowed us  to reduce the possibilities  which were initialized
in step~(5).

(9) It happens that the previous steps do  not rule out all but one tree. In
these  cases we  are always  left with  two or  four possible  trees and  it
turns  out that  the  remaining trees  are the  same  modulo permutation  of
algebraically conjugate characters  in the block. In such a  case any of the
possible  trees  is correct  for  some  choice  of  the modular  system  (or
equivalently,  for  some  choice  of  identification  of  certain  conjugacy
classes).

\subsection*{The cases to consider.}
The simple  group $\GG$ has  order $2^{36}\cdot 3^9\cdot 5^2\cdot  7^2 \cdot
11\cdot13 \cdot  17\cdot 19$, the orders  of the related groups  we consider
have additional factors $2$ or $3$. We handle the following cases where $G =
\GG$:
\begin{itemize}
\item $(2^2\times 3).G$ and $l = 5, 7, 11, 13, 17, 19$.
\item $2^2.G$ and $l = 3$.
\item $(2^2\times 3).G.2$ and $l = 5, 7, 11, 13, 17, 19$.
\item $2^2.G.2$ and $l = 3$.
\item $2^2.G.3$ and $l =  5, 7, 11, 13, 17, 19$.
\item $2^2.G.S_3$ and $l =  5, 7, 11, 13, 17, 19$.
\end{itemize}
Note  that in  a larger  extension we  also see  the blocks  of the  smaller
extensions which are quotient groups. The special cases for $l=3$ are needed
because  the blocks  with non-trivial  cyclic defect  have no  longer cyclic
defect  in  the  larger  extension.  For $l=2$  there  are  no  blocks  with
non-trivial cyclic defect in any  of these groups. In $2^2.G.3$, $2^2.G.S_3$
there are no blocks with non-trivial cyclic defect for $l=3$.

\subsection{The Brauer trees}

Here is  the list of  all Brauer trees  we obtain. We  give the name  of the
character tables as they can be  accessed in \GAP\ via the \CTblLib-package.
The nodes  of the  trees are  labeled by the  position of  the corresponding
character in that table. A notation $(i+j)$ denotes a sum of the exceptional
characters in a block.

All characters that do not appear here are either not in a block with cyclic
defect (this happens for $p \in \{3,5,7\}$) or of defect zero.

Some trees are only determined up to algebraic conjugacy, that is any of the
possible trees is correct with respect  to some choice of $p$-modular system
or to some  choice of labeling for certain conjugacy  classes. In such cases
we mention the  permutation group on the characters in  the block induced by
Galois automorphisms.

If the characters of a block are not faithful we indicate the order $|K|$ of
their kernel.

\tikzset{
  vertex/.style = {}
}

{\small


\subsubsection{\GAP\ table \texttt{(2\^{}2x3).2E6(2)}, $l = 5$}

\mbox{}

 \raisebox{0.43em}{$\quad (|K| = 4)$}

\subsubsection{\GAP\ table \texttt{(2\^{}2x3).2E6(2)}, $l = 7$}

\mbox{}

%
 \raisebox{0.43em}{modulo $\langle (787,788), (798,799) \rangle$} \raisebox{0.43em}{$\quad (|K| = 6)$}

\subsubsection{\GAP\ table \texttt{(2\^{}2x3).2E6(2)}, $l = 11$}

\mbox{}

%
 \raisebox{0.43em}{$\quad (|K| = 2)$}

\subsubsection{\GAP\ table \texttt{(2\^{}2x3).2E6(2)}, $l = 13$}

\mbox{}

%
\\
$\quad\quad$  \raisebox{0.43em}{modulo $\langle (878,880), (892,894) \rangle$} \raisebox{0.43em}{$\quad (|K| = 2)$}

\subsubsection{\GAP\ table \texttt{(2\^{}2x3).2E6(2)}, $l = 17$}

\mbox{}

%
\\
$\quad\quad$  \raisebox{0.43em}{modulo $\langle (838,840) \rangle$} \raisebox{0.43em}{$\quad (|K| = 2)$}

\subsubsection{\GAP\ table \texttt{(2\^{}2x3).2E6(2)}, $l = 19$}

\mbox{}

%
\\
$\quad\quad$  \raisebox{0.43em}{modulo $\langle (892,894) \rangle$} \raisebox{0.43em}{$\quad (|K| = 2)$}

\subsubsection{\GAP\ table \texttt{2\^{}2.2E6(2)}, $l = 3$}

\mbox{}

%
 \raisebox{0.43em}{$\quad (|K| = 4)$}

\subsubsection{\GAP\ table \texttt{(2\^{}2x3).2E6(2).2}, $l = 5$}

\mbox{}

%

\subsubsection{\GAP\ table \texttt{(2\^{}2x3).2E6(2).2}, $l = 7$}

\mbox{}

%
 \raisebox{0.43em}{modulo $\langle (539,540), (550,551) \rangle$}

\subsubsection{\GAP\ table \texttt{(2\^{}2x3).2E6(2).2}, $l = 11$}

\mbox{}

%

\subsubsection{\GAP\ table \texttt{(2\^{}2x3).2E6(2).2}, $l = 13$}

\mbox{}

%
\\
$\quad\quad$  \raisebox{0.43em}{modulo $\langle (630,632), (644,646) \rangle$}

\subsubsection{\GAP\ table \texttt{(2\^{}2x3).2E6(2).2}, $l = 17$}

\mbox{}

%
\\
$\quad\quad$  \raisebox{0.43em}{modulo $\langle (590,592) \rangle$}

\subsubsection{\GAP\ table \texttt{(2\^{}2x3).2E6(2).2}, $l = 19$}

\mbox{}

%
\\
$\quad\quad$  \raisebox{0.43em}{modulo $\langle (644,646) \rangle$}

\subsubsection{\GAP\ table \texttt{2\^{}2.2E6(2).2}, $l = 3$}

\mbox{}

%
 \raisebox{0.43em}{$\quad (|K| = 2)$}

\subsubsection{\GAP\ table \texttt{2\^{}2.2E6(2).3}, $l = 5$}

\mbox{}

%

\subsubsection{\GAP\ table \texttt{2\^{}2.2E6(2).3}, $l = 7$}

\mbox{}

%
 \raisebox{0.43em}{modulo $\langle (315,316), (326,327) \rangle$}

\subsubsection{\GAP\ table \texttt{2\^{}2.2E6(2).3}, $l = 11$}

\mbox{}

%

\subsubsection{\GAP\ table \texttt{2\^{}2.2E6(2).3}, $l = 13$}

\mbox{}

%
\\
$\quad\quad$  \raisebox{0.43em}{modulo $\langle (322,323) \rangle$}

\subsubsection{\GAP\ table \texttt{2\^{}2.2E6(2).3}, $l = 17$}

\mbox{}

%

\subsubsection{\GAP\ table \texttt{2\^{}2.2E6(2).3}, $l = 19$}

\mbox{}

%
\\
$\quad\quad$  \raisebox{0.43em}{modulo $\langle (322,323) \rangle$}

\subsubsection{\GAP\ table \texttt{2\^{}2.2E6(2).3.2}, $l = 5$}

\mbox{}

%

\subsubsection{\GAP\ table \texttt{2\^{}2.2E6(2).3.2}, $l = 7$}

\mbox{}

%

\subsubsection{\GAP\ table \texttt{2\^{}2.2E6(2).3.2}, $l = 11$}

\mbox{}

%

\subsubsection{\GAP\ table \texttt{2\^{}2.2E6(2).3.2}, $l = 13$}

\mbox{}

%
\\
$\quad\quad$  \raisebox{0.43em}{modulo $\langle (355,357) \rangle$}

\subsubsection{\GAP\ table \texttt{2\^{}2.2E6(2).3.2}, $l = 17$}

\mbox{}

%

\subsubsection{\GAP\ table \texttt{2\^{}2.2E6(2).3.2}, $l = 19$}

\mbox{}

%

}

\bibliographystyle{alpha}
\bibliography{twe62}

\end{document}